\magnification 1200

%
%
\newdimen\FigSize	\FigSize=.9\hsize 
%
\newskip\abovefigskip	\newskip\belowfigskip
\gdef\epsfig#1;#2;{\par\vskip\abovefigskip\penalty -500
   {\everypar={}\epsfxsize=#1\noindent
    \centerline{\epsfbox{#2}}}%
    \vskip\belowfigskip}%
%
\newskip\figtitleskip
\gdef\tepsfig#1;#2;#3{\par\vskip\abovefigskip\penalty -500
   {\everypar={}\epsfxsize=#1\noindent
    \vbox
      {\centerline{\epsfbox{#2}}\vskip\figtitleskip
       \centerline{\figtitlefont#3}}}%
    \vskip\belowfigskip}%
%
\newcount\FigNr	\global\FigNr=0
\gdef\nepsfig#1;#2;#3{\global\advance\FigNr by 1
   \tepsfig#1;#2;{Figure\space\the\FigNr.\space#3}}%
%
%
%
\gdef\ipsfig#1;#2;{
   \midinsert{\everypar={}\epsfxsize=#1\noindent
	      \centerline{\epsfbox{#2}}}%
   \endinsert}%
%
\gdef\tipsfig#1;#2;#3{\midinsert
   {\everypar={}\epsfxsize=#1\noindent
    \vbox{\centerline{\epsfbox{#2}}%
          \vskip\figtitleskip
          \centerline{\figtitlefont#3}}}\endinsert}%
%
\gdef\nipsfig#1;#2;#3{\global\advance\FigNr by1%
  \tipsfig#1;#2;{Figure\space\the\FigNr.\space#3}}%
\newread\epsffilein    
\newif\ifepsffileok    
\newif\ifepsfbbfound   
\newif\ifepsfverbose   
\newdimen\epsfxsize    
\newdimen\epsfysize    
\newdimen\epsftsize    
\newdimen\epsfrsize    
\newdimen\epsftmp      
\newdimen\pspoints     
\pspoints=1bp          
\epsfxsize=0pt         
\epsfysize=0pt         
\def\epsfbox#1{\global\def\epsfllx{72}\global\def\epsflly{72}%
   \global\def\epsfurx{540}\global\def\epsfury{720}%
   \def\lbracket{[}\def\testit{#1}\ifx\testit\lbracket
   \let\next=\epsfgetlitbb\else\let\next=\epsfnormal\fi\next{#1}}%
\def\epsfgetlitbb#1#2 #3 #4 #5]#6{\epsfgrab #2 #3 #4 #5 .\\%
   \epsfsetgraph{#6}}%
\def\epsfnormal#1{\epsfgetbb{#1}\epsfsetgraph{#1}}%
\def\epsfgetbb#1{%
%
%
\openin\epsffilein=#1
\ifeof\epsffilein\errmessage{I couldn't open #1, will ignore it}\else
%
%
   {\epsffileoktrue \chardef\other=12
    \def\do##1{\catcode`##1=\other}\dospecials \catcode`\ =10
    \loop
       \read\epsffilein to \epsffileline
       \ifeof\epsffilein\epsffileokfalse\else
%
%
          \expandafter\epsfaux\epsffileline:. \\%
       \fi
   \ifepsffileok\repeat
   \ifepsfbbfound\else
    \ifepsfverbose\message{No bounding box comment in #1; using defaults}\fi\fi
   }\closein\epsffilein\fi}%
%
%
\def\epsfsetgraph#1{%
   \epsfrsize=\epsfury\pspoints
   \advance\epsfrsize by-\epsflly\pspoints
   \epsftsize=\epsfurx\pspoints
   \advance\epsftsize by-\epsfllx\pspoints
%
%
   \epsfxsize\epsfsize\epsftsize\epsfrsize
   \ifnum\epsfxsize=0 \ifnum\epsfysize=0
      \epsfxsize=\epsftsize \epsfysize=\epsfrsize
%
%
     \else\epsftmp=\epsftsize \divide\epsftmp\epsfrsize
       \epsfxsize=\epsfysize \multiply\epsfxsize\epsftmp
       \multiply\epsftmp\epsfrsize \advance\epsftsize-\epsftmp
       \epsftmp=\epsfysize
       \loop \advance\epsftsize\epsftsize \divide\epsftmp 2
       \ifnum\epsftmp>0
          \ifnum\epsftsize<\epsfrsize\else
             \advance\epsftsize-\epsfrsize \advance\epsfxsize\epsftmp \fi
       \repeat
     \fi
   \else\epsftmp=\epsfrsize \divide\epsftmp\epsftsize
     \epsfysize=\epsfxsize \multiply\epsfysize\epsftmp   
     \multiply\epsftmp\epsftsize \advance\epsfrsize-\epsftmp
     \epsftmp=\epsfxsize
     \loop \advance\epsfrsize\epsfrsize \divide\epsftmp 2
     \ifnum\epsftmp>0
        \ifnum\epsfrsize<\epsftsize\else
           \advance\epsfrsize-\epsftsize \advance\epsfysize\epsftmp \fi
     \repeat     
   \fi
%
%
   \ifepsfverbose\message{#1: width=\the\epsfxsize, height=\the\epsfysize}\fi
   \epsftmp=10\epsfxsize \divide\epsftmp\pspoints
   \vbox to\epsfysize{\vfil\hbox to\epsfxsize{%
      \includegraphics{#1}%
      \hfil}}%
\epsfxsize=0pt\epsfysize=0pt}%
%
%
{\catcode`\%=12 \global\let\epsfpercent=
%
%
\long\def\epsfaux#1#2:#3\\{\ifx#1\epsfpercent
   \def\testit{#2}\ifx\testit\epsfbblit
      \epsfgrab #3 . . . \\%
      \epsffileokfalse
      \global\epsfbbfoundtrue
   \fi\else\ifx#1\par\else\epsffileokfalse\fi\fi}%
%
%
\def\epsfgrab #1 #2 #3 #4 #5\\{%
   \global\def\epsfllx{#1}\ifx\epsfllx\empty
      \epsfgrab #2 #3 #4 #5 .\\\else
   \global\def\epsflly{#2}%
   \global\def\epsfurx{#3}\global\def\epsfury{#4}\fi}%
%
%
\def\epsfsize#1#2{\epsfxsize}%
%
%

\epsfverbosetrue			
\abovefigskip=\baselineskip		
\belowfigskip=\baselineskip		
\global\let\figtitlefont\bf		
\global\figtitleskip=.5\baselineskip	

\font\tenmsb=msbm10   
\font\sevenmsb=msbm7
\font\fivemsb=msbm5
\newfam\msbfam
\textfont\msbfam=\tenmsb
\scriptfont\msbfam=\sevenmsb
\scriptscriptfont\msbfam=\fivemsb
\def\Bbb#1{\fam\msbfam\relax#1}
\let\nd\noindent 
\def\qed{\hbox{\hskip 6pt\vrule width6pt height7pt depth1pt \hskip1pt}}
\def\natural{{\rm I\kern-.18em N}}

\def\integer{{\rm Z\kern-.32em Z}}
\def\chix{{\raise.5ex\hbox{$\chi$}}}
\def\Z{{\Bbb Z}}
\def\real{{\rm I\kern-.2em R}}
\def\R{{\Bbb R}}

\def\complex{\kern.1em{\raise.47ex\hbox{
            $\scriptscriptstyle |$}}\kern-.40em{\rm C}}

\def\A{{\cal A}}
\def\O{{\cal O}}
\def\G{{\hat {\cal O}}}

\def\a{\alpha}
\def\b{\beta}
\def\vs#1 {\vskip#1truein}
\def\hs#1 {\hskip#1truein}
  \hsize=6truein        \hoffset=.25truein 
  \vsize=8.8truein      
  \pageno=1     \baselineskip=12pt
  \parskip=5 pt         \parindent=20pt
  \overfullrule=0pt     \lineskip=0pt   \lineskiplimit=0pt
  \hbadness=10000 \vbadness=10000 
\pageno=0

\footline{\ifnum\pageno=0\hss\else\hss\tenrm\folio\hss\fi}
\hbox{}
\vskip 1truein\centerline{{\bf AN ALGEBRAIC INVARIANT FOR SUBSTITUTION 
TILING SYSTEMS}}
\vskip .5truein\centerline{by}
\centerline{Charles Radin${}^1$ and Lorenzo Sadun${}^2$} 
\footnote{}{1\ Research supported in part by NSF Grant No. DMS-9531584\hfil}
\footnote{}{2\ Research supported in part by NSF Grant No. DMS-9626698\hfil}
\vskip .2truein\centerline{Mathematics Department}
\centerline{University of Texas}
\centerline{Austin, TX\ \ 78712}
\vs.1
\centerline{radin@math.utexas.edu and sadun@math.utexas.edu}
\vs.5
\centerline{{\bf Abstract}}
\vs.1 \nd
We consider tilings of Euclidean spaces by polygons or polyhedra, in
particular, tilings made by a substitution process, such as the
Penrose tilings of the plane. We define an isomorphism invariant
related to a subgroup of rotations and compute it for various
examples. We also extend our analysis to more general dynamical
systems.
\vs.8 \nd
Running head: INVARIANT FOR SUBSTITUTION TILING SYSTEMS
\vs.1 \nd
1991 AMS Classification: 52C22, 52C20
\vfill\eject 
\nd
{\bf 1.\ Introduction, Definitions and Statement of Results} 
\vs.1 
This paper concerns tilings of Euclidean
spaces by polygons or polyhedra, more specifically, tilings made by a
``substitution process''.  Given a substitution rule, the set of resultant
tilings is a topological space with an action of the Euclidean group,
hence a dynamical system.  We
develop here an algebraic invariant that helps determine when two tiling
systems are equivalent {\it as dynamical systems}.  In this introduction
we define the notions of ``substitution tiling system'' and of equivalence
between two such systems, and state what the invariant is.  In subsequent
sections we analyze the invariant, in particular we show its
use in distinguishing between substitution tilings. 

Our eventual goal is to associate certain groups to
substitution tilings of Euclidean $m$-space. These groups, subgroups of
$SO(m)$, are generated by the relative orientations of tiles in the
tilings, depend on the specific tiling $x$ and on some specific
choices (indexed by an integer $j$), and are denoted $\O_j(x)$.  
Although $\O_j(x)$ depends on $j$ and $x$, the dependence is quite controlled.
If $x$ and $x'$ are different tilings with the same substitution rule we
will show that, under some mild hypotheses, $\O_j(x)$ and $\O_{j'}(x')$ 
are conjugate as subgroups of $SO(m)$. Even without the mild hypotheses, 
they are conjugate (in $SO(m)$)
to subgroups of one another, a condition we call ``c-equivalence''.  
So we can associate to a substitution tiling system the common conjugacy 
class (or c-equivalence class) of the groups associated to the tilings
in the system. 

What will remain, then, is to show that this conjugacy (or
c-equivalence) class can be considered an invariant in a natural
sense.  That is, we will show that two substitution tiling systems
that are equivalent as dynamical systems have the same class of
groups.  We will do this by finding a dynamical description of the
class.  For each $\epsilon>0$ and each tiling $x$ we will define a
group $\O^\epsilon(x)$ using dynamical information only.  For
$\epsilon$ sufficiently small, and for almost every $x$, we show that
$\O^\epsilon(x)$ is conjugate to (or c-equivalent to) $\O_j(x)$ for
some, and hence all, choices $j$.  The class of $\O^\epsilon(x)$ is
thus the same as the class of $\O_j(x)$.  Since $\O^\epsilon(x)$ is
defined using data that is preserved by dynamical equivalence, the
class of $\O^\epsilon(x)$ is a dynamical invariant.

Note that the group $\O_j(x)$ depends only on the geometry of the
tiling $x$.  Since the class of $\O_j(x)$ is the same for every tiling
$x$ with the given substitution rule,
we can obtain information about a substitution tiling
system by looking at any single tiling in it.  So if two substitution tilings $x$
and $x'$ give rise to groups $\O_j(x)$ and $\O_{j'}(x')$ that are not
conjugate (or c-equivalent), then $x$ and $x'$ cannot belong to
equivalent substitution tiling systems.

Before defining substitution tiling systems in general, we present
an example.  Hopefully, the general definitions will be clearer with
this example in mind.  The  ``pinwheel'' tiling
of the plane [Ra1] is made as follows. Consider the triangles of
Fig.\ 1. Divide one of them  into five
small triangles as in Fig.\ 2 and expand the figure about the origin by a
linear factor of $\sqrt 5$, producing 5 triangles congruent to the originals.
\vs1
\epsfig .5\hsize; 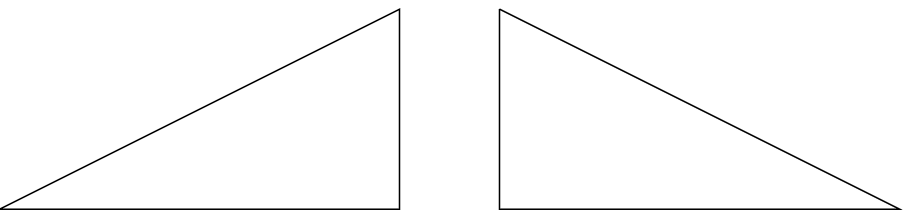 ; 
\vs.1
\centerline{Figure 1. Two ``pinwheel'' tiles}
\vs.1
\vbox{\epsfig .5\hsize; 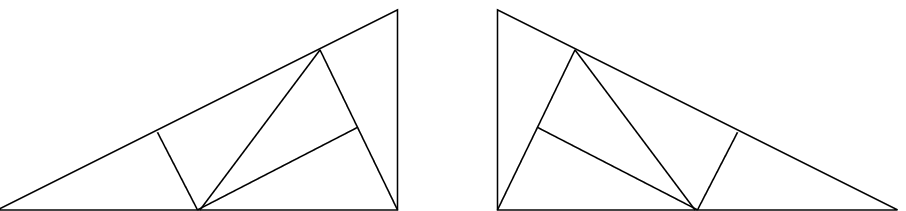 ;}
\vs.1
\centerline{Figure 2. The substitution for pinwheel tilings}
\vs-2.6
\hs-.2 \vbox{\epsfig 1\hsize; 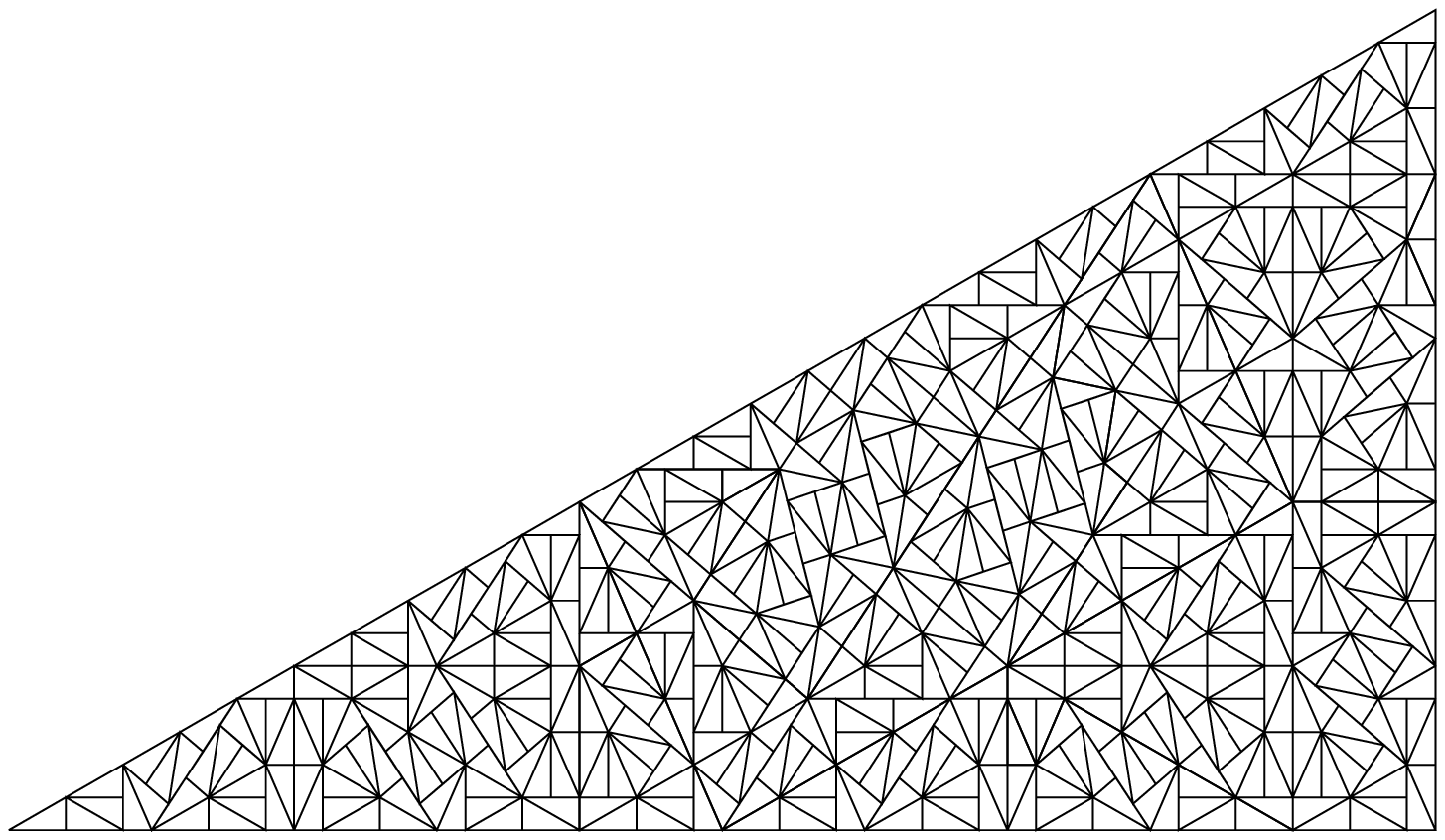 ;}
\hs.2
\centerline{Figure 3. Part of a pinwheel tiling}

Repeat this two-step
procedure $\phi$ simultaneously for all the triangles of the figure, then
again, an infinite number of times, producing a (pinwheel) tiling $C$ of
the plane, a portion of which appears in Fig.\ 3. Such tilings have a
hierarchical structure which is of interest for various reasons; in
particular it leads to interesting behavior of the relative
orientations of tiles within a tiling [Ra3]. For background on related
recent work see [AnP, CEP, DwS, G-S, Kel, Ken, LaW, Min, Moz, Ra3, Rob,
Sad, Sch, Sen, Sol, Tha] and references therein.

\vs.05 \nd
{\bf Substitution Tiling Systems}
\vs.05
With the pinwheel example in mind, we now address substitution tiling
systems in general.
Let $\A$ be a nonempty finite collection of polyhedra in $m$
(typically 2 or 3) dimensions. Let $X(\A)$ be the set of all tilings
of Euclidean space by congruent copies, which we will call tiles, of
the elements of (the ``alphabet'') $\A$. We label the 
``types'' of tiles by the elements
of $\A$.  We endow $ X(\A)$ with the metric

$$d(x, y)\equiv \sup_{n}{1\over n}m_{H}[B_{n}(\partial x),B_{n}(\partial
x')], \eqno 1)$$

\nd where ${B}_{n}(\partial x)$ denotes the intersection of two sets:
the closed ball ${B}_n$ of radius $n$ centered at the origin of the
Euclidean space and the union $\partial x$ of the boundaries $\partial
a$ of all tiles $a$ in $x$. $m_{H}$ is the Hausdorff metric on compact
sets defined as follows. Given two compact subsets $A$ and $B$ of $\R^m$,
$m_{H}[A,B] = \max
\{ {\tilde d} (A,B), {\tilde d} (B,A)\}$, where
$${\tilde d} (A,B) =  \sup_{a \in A}\inf_{b\in B} ||a - b||, \eqno 2)$$

\nd with $||w||$ denoting the usual Euclidean norm
of $w$.  Although the metric $d$ depends on the location of the
origin, the topology induced by $d$ is translation invariant.  A
sequence of tilings converges in the metric $d$ if and only if its
restriction to every compact subset of $\R^m$ converges in $m_H$.
It is not hard to show [RaW] that $X(\A)$ (which is automatically
nonempty in our applications) is compact and that the natural action of the connected
Euclidean group $G_E$ on $ X(\A)$, $(g,x)\in G_E\times
X(\A)\longrightarrow T^gx\in X(\A)$, is continuous. 

A ``substitution tiling system'' is a closed subset $X_\phi \subset
X(\A)$ satisfying some additional conditions. 
To understand these conditions we first need the notion of
``patches''.  A patch is a (finite or infinite) subset of an element
$x\in X(\A)$; the set of all patches for a given alphabet will be
denoted by $W$. Next we need, as for the pinwheels, an auxiliary
``substitution function'' $\phi$, a map from $W$ to $W$, with the
following properties:
\vs.1
\vs0 \hs.2 i) \phantom{ii}There is some constant $|\phi|>1$ such 
that, for any $g\in G_E$ and $x\in X$, 
\vskip -6pt \hs.2 \phantom{iii) }$\phi[T^gx]=T^{\phi(g)}x$, where 
$\phi(g)$ is the conjugate of $g$ by the similarity of
\vskip -6pt \hs.2 \phantom{iii) }Euclidean space consisting of
stretching about the origin by $|\phi|$.
\vskip -6pt \hs.2 ii) \phantom{i}For each tile $a\in\A$ and for each 
$n \ge 1$, the union of the tiles in $\phi^n a$ is
\vskip -6pt \hs.2 \phantom{iii) }congruent to $|\phi|^n a$, and these
tiles meet full face to full face.
\vskip -6pt \hs.2 iii) For each tile $a\in\A$, $\phi a$ contains at
least one tile of each type.
\vs.1 
Condition ii) is quite strong. It is satisfied by the pinwheel tilings 
only if we add additional vertices at midpoints of the legs of
length 2, creating boundaries of 4 edges.  A similar (minor) adjustment is
needed for other examples in this paper.  Even
with such adjustments however, condition ii) is not satisfied by the kite \&
dart tilings [Gar], or those which mimic substitution tilings using
so-called edge markings [G-S, Moz, Ra3]. It is to handle such examples 
that we introduce the general development of \S 3.
\vs.1 \nd
{\bf Definition.} For a given alphabet $\A$ of polyhedra and
substitution function $\phi$  the ``substitution tiling system''
is the pair $\{X_\phi, T\}$, where 
$X_\phi \subset X(\A)$ is the compact subset of
those tilings $x$ with the property that every finite
subpatch of $x$ is congruent to a subpatch of $\phi^n a$ for some
$n>0$ and some $a\in \A$, and $T$ is the natural action of
$G_E$ on $X_\phi$. (For simplicity we often refer to
$X_\phi$ as a substitution tiling system.)
\vs.1 

One
planar example of a substitution tiling system 
is based on the pinwheel substitution of Fig.~2.
A slight variant of the pinwheel is defined by the 1-3-$\sqrt{10}$ right
triangle and its reflection, and the substitution of Fig.\ 4.
\vs.1
\vbox{\epsfig .8\hsize; 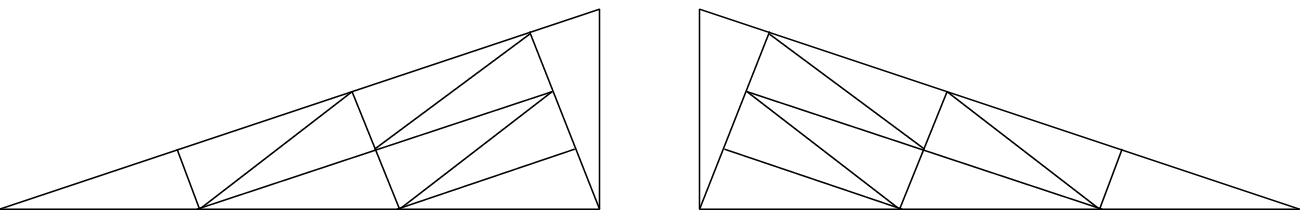 ;}
\vs.1
\centerline{Figure 4. The substitution for pinwheel variant tilings}
\vs.1 \nd

Two further special conditions which we will occasionally impose are:
\vs.1 \nd
iv) A tiling in $X_\phi$ can only be tiled in
one way by supertiles of level $n$, for any 
\vs0 \nd
\phantom{iv) }$n\ge 1$.
\vs.05
\nd v) For every $a \in \A$, there exists $n_a>0$ such that 
$\phi^{n_a}a$ contains a tile $a'$, of the 
\phantom{v) }same type as $a$, and parallel to $a$.
\vs.1
We note here that with the convention that
patches of the form $\phi^n a$ are called ``supertiles'' of ``level''
$n$ and ``type'' $a$, it is easy to show by a diagonal argument that,
for each $n\ge 0$, each tiling $x$ is tiled by supertiles of level
$n$ [Ra3]. A supertile of level 4 for the pinwheel is shown in Fig.\ 3.

Finally, let
\vs.05
$S^\epsilon(x)= \{y\in X_{\phi}\,:\,d(\phi^ny,\phi^nx)<\epsilon
\hbox{ for all }n\ge 0,
\hbox{ and }d(\phi^ny,\phi^nx)\mathop{\longrightarrow}
\limits_{n\to\infty}0\}.$
\vs.05 \nd
We call such a family of sets a ``local contracting direction (at $x$)''.

Our goal is to define a notion of equivalence for substitution tiling
systems, and an invariant for that equivalence. For the equivalence we
use:
\vs.1 \nd
{\bf Definition.} The substitution tiling systems $(X_{\phi^1},
T^1)$ and $(X_{\phi^2}, T^2)$ are ``equivalent'' if there are
subsets $Y_j\subset X_{\phi^j}$, invariant under $T^j$ and of
measure zero with respect to all translation invariant Borel
probability measures on $X_{\phi^j}$, 
and a one-to-one, onto,
Borel bimeasureable map $\tau\,:\,X_{\phi^1}\!-\!Y_1\to X_{\phi^2}\!-\!Y_2$,
such that: \vs0 \hs.2
a) $\tau\circ T^1=T^2\circ \tau$; 
\vs0 \hs.2
b) for each $x\in X^1\!-\!Y_1$, $\epsilon
>0$ and $\epsilon' >0$, there exist $\tilde \epsilon >0$ and $\tilde
\epsilon' >0$ such 
\vskip -6pt \hs.2 \phantom{b) }that $\tau[S^{\tilde \epsilon}(x)]
\subset S^{\epsilon}(\tau x)$, and $\tau^{-1}[S^{\tilde
\epsilon'}(\tau x)] \subset S^{\epsilon'}(x)$. \vs0 
\nd We call such a map $\tau$
an ``isomorphism''.

This notion of equivalence is stronger than simply intertwining the
actions of $G_E$. This is appropriate; it has been known at least
since [CoK] that substitution subshifts show almost none of their
richness if considered merely as subshifts. So in classifying tilings
that have a hierarchical structure we make some feature of that
hierarchical structure part of our notion of equivalence. 

To define an invariant we extract information from the local
contracting directions.  Since the local contracting directions are
preserved by equivalence, such information is manifestly invariant. We
define here the invariant.  In later sections we relate it to directly
computable quantities (the $\O_j(x)$) and demonstrate its use in
distinguishing between tiling dynamical systems.
 
Consider $G_E$ as the semidirect product of $SO(m)$ with $\R^m$, 
with $g=(r,t)$ denoting a
rotation $r$ about the origin followed by a translation $t$. 
Then consider, for any 
substitution tiling system $X_{\phi}$ and $\epsilon >0$:
\vs.1
$R^\epsilon(x) =
\{r\in SO(m)\,:\,\hbox{there exists }t\hbox{ such that }
T^{(r,t)}x\in S^\epsilon(x)\}$ 
\vs.1 \nd
Now let $\O^\epsilon(x)$ be the subgroup of $SO(m)$ generated by
$R^\epsilon(x)$.  
The corollary to Theorem 2 shows that the conjugacy class of
$\O^\epsilon(x)$ is independent of $x$ and $\epsilon$ (when small
enough) for substitution tiling
systems satisfying iv) and v). The conjugacy class of $\O^\epsilon(x)$
is therefore an invariant of the tiling dynamical system, not just
a feature of the individual tiling $x$.

\nd{\bf 2.\ The group of relative orientations}

The group $\O^\epsilon(x)$ generated by $R^\epsilon(x)$ is not
directly computable.  In this section we remedy this by constructing,
for a substitution tiling system, a more easily computable group
$\O_j(x)$ related to the relative orientations of the tiles in the
single tiling $x$. The group $\O^\epsilon(x)$ is then shown to be
conjugate to $\O_j(x)$.

Given a tiling $x$ and some tile $a$ of type $j$ in it, let
$R_j(a,x)\subset SO(m)$ be the set of relative orientations with
respect to $a$ of the tiles of type $j$ in $x$; that is, $R_j(a,x)$ is
the set of rotations of $x$ which bring a tile of type $j$ parallel to
(the fixed) $a$. The group generated by $R_j(a,x)$ is easily seen to
be generated by the relative orientations between {\it all} pairs of tiles
of type $j$ in $x$; in particular it is independent of $a$, and we
denote it by $\O_j(x)$. Furthermore,
\vs.1 \nd 
{\bf Lemma 1.}\ If $x'$ has a tile $a'$ of type $j$ parallel to the
tile $a$ in $x$, then $\O_j(x)=\O_j(x')$. For any $\tilde x$,
$\O_j(\tilde x)$ is conjugate to $\O_j(x)$.
\vs.1 \nd
Proof.\ First note that $\O_j(x)$ is generated by the relative
orientations between $a$ and all other tiles of type $j$ in $x$. So
consider the generator $g$ of $\O_j(x)$ which is the relative
orientation of a tile $c$ with respect to $a$ in $x$.  We will show
that $g \in \O_j(x')$, from which it follows that $\O_j(x) \subset
\O_j(x')$.  By symmetry, we would then have $\O_j(x')
\subset \O_j(x)$, and hence $\O_j(x)=\O_j(x')$.

 From the definition of substitution tilings, the tiles $a$ and $c$
can be thought of as belonging to some supertile $A$ of level $n$
(although not all of $A$ need exist in $x$). Since $x'$ is tiled by
supertiles of level $n$, there is a supertile $A'$ of level $n$ in $x'$
containing a pair of tiles, $a''$ and $c''$, which have the same
positions relative to $A'$ as do $a$ and $c$ relative to $A$. See
Fig.~5.
\vs.1
\epsfig 1\hsize; 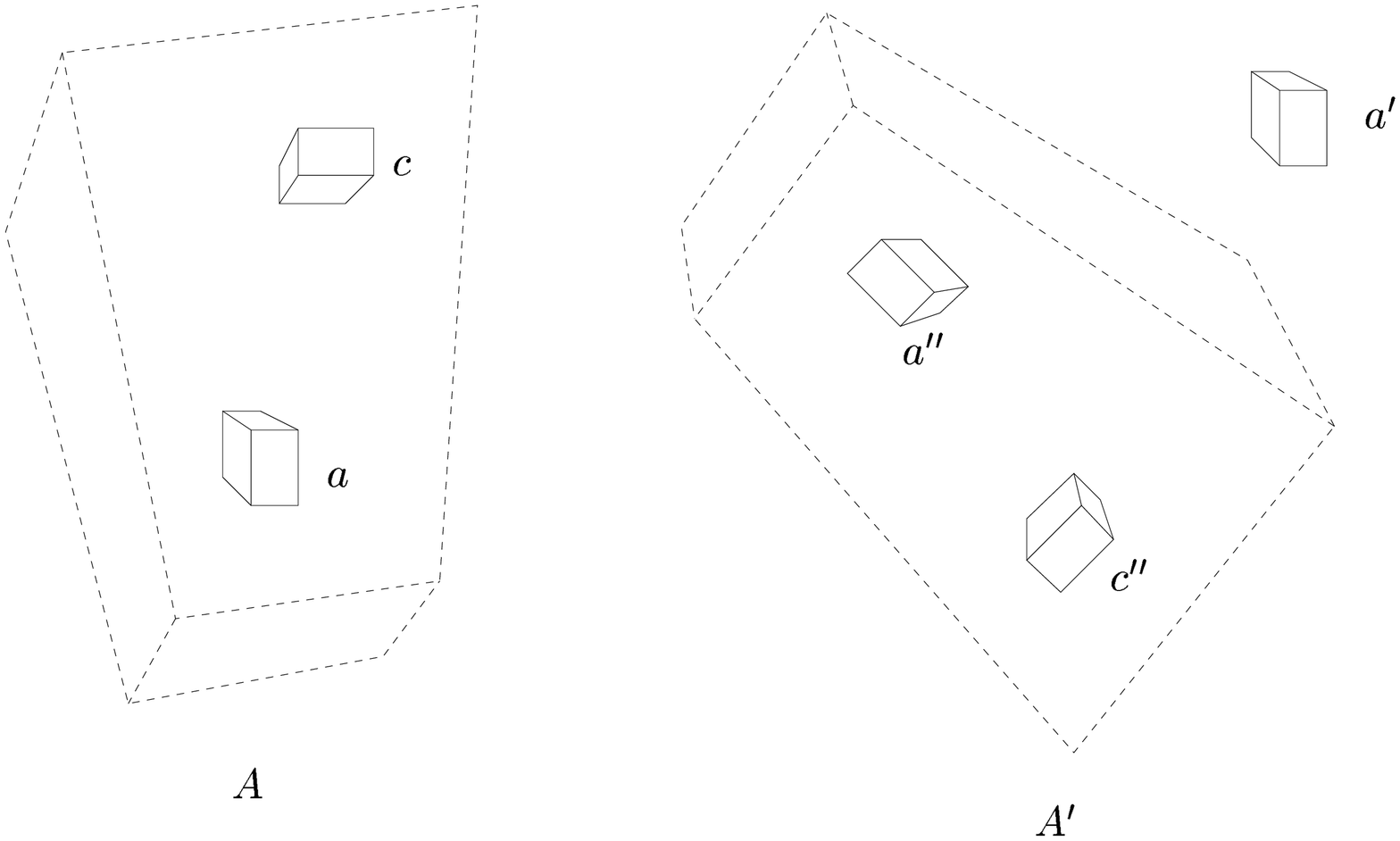 ; 
\vs.1
\hs.5 part of tiling $x$\hs2 part of tiling $x'$
\vs.1
\centerline{Figure 5}
\vs.1
Let $g'$ be the relative
orientation of $c''$ with respect to $a''$. Then $g=R^{-1}g'R$ where
$R$ is the relative orientation of $A'$ with respect to $A$. But $R$
is then also the relative orientation of $a''$ with respect to $a$,
which is the same as that of $a''$ with respect to $a'$, and $R$ is
thus an element of $\O_j(x')$. But then $g$ is an element of $\O_j(x')$,
as claimed. 

If $\tilde x$ is any tiling at all, $\O_j(x)$ and $\O_j(\tilde x)$
are conjugate by an element of $SO(m)$, namely a rotation which 
makes a tile of type $j$ in $\tilde x$ parallel to one in $x$. \qed

Finally we consider the dependence of $\O_j(x)$ on $j$. 
\vs.1 \nd
{\bf Definition.} Two subgroups of $SO(m)$ are ``c-equivalent'' if each is
conjugate (in $SO(m)$) to a subgroup of the other. (Note that in
$SO(2)$ c-equivalence is the same as identity.)
\vs.1
\nd {\bf Lemma 2.}\ For any tilings $x, \tilde x \in X_\phi$ and tile types 
$j$ and $k$, $\O_j(x)$ is c-equivalent to $\O_k(\tilde x)$.
\vs.1 \nd
Proof. By Lemma 1 it is sufficient to show that $\O_j(x)$ and $\O_k(x)$
are c-equivalent.  Consider any two tiles $a$ and $b$ of type $j$ in $x$, and let
$g$ be the relative orientation of $b$ with respect to $a$.
After one substitution $a$ and $b$ give rise to
tiles $a'$ and $b'$ of type $k$ in the tiling $\phi x$. The relative
orientation of $b'$ with respect to $a'$ is again $g$, since $g$ takes
each part of $b$ onto the corresponding part of $a$.  Applying this
construction to all the generators of $\O_j(x)$, we see that $\O_j(x)$
is a subgroup of $\O_k(\phi x)$. Similarly, $\O_k(x)$ is a subgroup of
$\O_j(\phi x)$. But $\O_j(x)$ and $\O_k(x)$ are conjugate to
$\O_j(\phi x)$ and $\O_k(\phi x)$, respectively, so $\O_j(x)$ and
$\O_k(x)$ are conjugate to subgroups of each other.  \qed
\vs.1
\nd {\bf Lemma 3.} Assume $\phi$ satisfies v).  Then $\O_j(x)
=\O_k(x)$.

\nd Proof. Let $n = \Pi_a n_a$. Since the tilings defined by $\phi^n$ are 
the same as those defined by $\phi$ we can, without loss of
generality, assume $n=1$, so that $\phi x$ contains tiles parallel to
every tile of $x$.  Then, by Lemma 1, $\O_j(x)=\O_j(\phi x)$ and $\O_k(x)
= \O_k(\phi x)$.  But we have shown that $\O_k(x) \subset \O_j(\phi x)$
and $\O_j(x)
\subset \O_k(\phi x)$, so $\O_j(x)=\O_k(x)$. \qed

To summarize: From lemmas 1 and 2 we can associate a subgroup of
$SO(m)$ to any substitution tiling system, uniquely defined up to
c-equivalence. If the substitution tiling system satisfies
v), lemma 3 shows that the group is uniquely defined up to
conjugacy.

Before we can use these groups as an invariant for equivalence of
substitution tiling systems we must refer to the relative orientations
in a more fundamental way. Our next goal is to connect this group with
the invariant introduced at the end of \S 1. The essential observation
is that, if tilings $x\ne y$ agree in some neighborhood of the origin
in Euclidean space, then $\phi x$ and $\phi y$ will agree in a larger
neighborhood of the origin, so we typically expect $d(\phi x,\phi y)<
d(x,y)$.  We are thus led to a quantity introduced earlier. For each
$x$ in the substitution tiling system $X_\phi$ and for each $\epsilon
> 0$, consider:
$$
S^\epsilon(x)=\{y\in X_{\phi}\,:\,d(\phi^ny,\phi^nx)<\epsilon
\hbox{ for all }n\ge 0,
\hbox{ and }d(\phi^ny,\phi^nx)\mathop{\longrightarrow}
\limits_{n\to\infty}0\}.
\eqno 3)
$$

\nd {\bf Theorem 1.} Assume a substitution tiling system $X_\phi$.
\vs0 
\item {a)} Given any $\epsilon>0$ there
exists $N>0$ such that ${B}_N(\partial y)={B}_N(\partial x)$
implies $y\in S^\epsilon(x)$.
\item{b)} There exists $\epsilon>0$ such that, for every $x\in X_\phi$,
$y\in S^\epsilon(x)$ and every tile $a \in x$ that meets the origin,
there is a tile $a' \in y$ that exactly coincides with $a$.
\vs.1
\nd Proof.\  a) is immediate from the form of the metric.  The proof of b)
requires the following two lemmas.

\nd {\bf Lemma 4.}  For every $N>0$ and every neighborhood $U$ of the
identity in $G_E$ there exists $\epsilon>0$ with the following
property: Let $x,x'\in X(\A)$ be any two tilings with
$d(x,x')<\epsilon$, and let $a$ be a tile of $x$ that is contained in
$B_N$. Then $x'$ contains a tile $a'$ of the form $T^ga$ where $g\in
U$.

\nd Proof.\ Let $z>0$ be such that for each $b\in \A$ some ball of
diameter $z$ lies in the interior of $b$. Fix some $\delta \in (0,z/3)$ and
define the heart $h_\delta(b)$ of $b\in \A$ as $\{p\in b\,:\, ||p-q||>
\delta\hbox{ for all }q\in \partial b\}$. By the corridor
$C_\delta(x)$ of a tiling $x$ we mean the complement of $\cup_{b\in
x}h_\delta(b)$.  Let $D$ be the largest of the diameters of all $b\in
\A$.  Without loss of generality, we can assume $N>D$.

With this notation we note that if $d(x,x')< \delta/(N+D)$ we have
$B_N(\partial x)\subset C_\delta(x')$ and $B_N(\partial x') \subset
C_\delta(x)$. So if $\delta \sim 0$ each tile in $x$ in $B_N$ is
closely approximated by some tile in $x'$ and vice versa. In
particular it now follows that for small enough $\epsilon$, if
$d(x,x')<\epsilon$ then the tiles $a'\in x'$ must be of the same type
as the tiles $a\in x$ they approximate, and in fact satisfy $a'=T^ga$
with $g\in U$. \qed

\nd {\bf Lemma 5.}  For each $N>0$ there is a
neighborhood $U_N$ of the identity in $G_E$ with the following
property:  If $a$ is a tile in $B_N$ and $g_1$ and $g_2$
are distinct elements of $U_N$,
then $T^{g_1}a$ and $T^{g_2}a$ overlap but are distinct. In
particular, it is impossible for $T^{g_1}a$ and $T^{g_2}a$ to both be
tiles in the same tiling.

\nd Proof.  This follows from the continuity of the action of $G_E$ on
tiles, and the fact that polyhedra do not admit infinitesimal
symmetries. \qed

We now return to the proof of Theorem 1. Pick
$N>D$ and let $U_N$
be as in Lemma 5. Pick a smaller bounded neighborhood $U \subset U_N$
of the identity of $G_E$ with the property that $(r,t) \in U$ implies
$(r,|\phi|t) \in U_N$.  Then pick $\epsilon$ small enough that Lemma 4
applies.

Let $a$ be a tile of $x$ containing the origin. By Lemma 4 there is a
tile $a'$ in $y$ of the form $T^g a $ with $g=(r,t) \in U\subset G_E$.
We will show that $t \ne 0$ implies that, for some $n$, $d(\phi^n x,
\phi^n y) > \epsilon$, while $t=0,\ r \ne 0$ implies that $\lim_{n \to
\infty} d(\phi^n x,\phi^ny) \ne 0$. This will complete the proof.

Note that $\phi^n a' = T^{(r,|\phi|^nt)}\phi^n a$.  If $t \ne 0$, pick
$n$ such that $(r,|\phi|^n t)$ is outside the neighborhood $U$ but
in $U_N$.   Let $\tilde a$ be
a tile of $\phi^n a$ containing the origin.  Then there is a tile
$\tilde a' = T^{(r,|\phi|^nt)}\tilde a$ in $\phi^n y$.  But by Lemma 5
this means there cannot be a tile of the form $T^g \tilde a$ in $y$
with $g \in U$.  By Lemma 4 this means that $d(\phi^nx,\phi^ny) \ge
\epsilon$. 

If $t=0$ then $\phi^n a' = T^{(r,0)} \phi^n a$. If $r\ne 0$, for every
tile $\tilde a\in \phi^na$ containing the origin there is a tile
$\tilde a'\in \phi^na'$ overlapping it and with relative orientation
$r$, which implies that the distance between $\phi^n x$ and $\phi^n y$
will not go to zero. \qed

Recall the following quantity from \S 1:
$$R^\epsilon(x) =
\{r\in SO(m)\,:\,\hbox{there exists }t\hbox{ such that }
T^{(r,t)}x\in S^\epsilon(x)\}.\eqno 4)$$ 
Let $\O^\epsilon(x)$ be the group
generated by $R^\epsilon(x)$. Assuming $\epsilon$ small enough for
Theorem 1b, we see that every $r \in R^\epsilon(x)$ is the
relative orientation of a tile of $x$ with respect to a corresponding
tile of $x$ near the origin.  By Theorem 1a, if $C$ is a region of
$x$ containing $B_N$, and if $C'$ is any region of $x$ congruent to $C$,
then $R^\epsilon(x)$ includes the relative orientation of $C'$ to
$C$.

Consider the following property.
\vs.1
\nd {\bf Property F.}\ The subset of tilings $x$, for which every 
fixed finite ball $B$ of Euclidean space is contained in some
supertile of finite level in $x$, is of full measure for every
translation invariant measure on $X_{\phi}$.
\vs.1 \nd
We will prove that Property F holds for a large class of interesting
systems, at least those satisfying condition iv).
This assumption, which implies that $\phi$ is a homeomorphism on
$X_{\phi}$, is satisfied by all known nonperiodic examples. In fact it
is automatically true for a system that contains
nonperiodic tilings and in which the tiles only appear in finitely
many orientations in any tiling [Sol].

If a tiling contains two or more regions each tiled by supertiles of level
$n$ for all $n\ge 0$, we call these regions supertiles of infinite level.
Recall that any tiling is tiled by supertiles of any finite level $n$. If a
ball in a tiling $x$ fails to lie in {\it any} supertile of any level $n$,
then $x$ is tiled by two or more supertiles of infinite level, with the
offending ball straddling a boundary.  (One can construct a pinwheel
tiling with two supertiles of infinite level as follows. Consider the
rectangle consisting of two supertiles of level $n-1$ in the middle of a
supertile of level $n$. For each $n\ge 1$ orient such a rectangle with its
center at the origin and its diagonal on the $x$-axis, and fill out
the rest of a (non-pinwheel) tiling $x_n$ by periodic extension. By
compactness this sequence has a convergent subsequence, which will be
a pinwheel tiling and which will consist of two supertiles of infinite
level.)

We now use the above to prove:
\vs.1 \nd
{\bf Lemma 6.}\ \ For a substitution tiling system satisfying iv), 
let $S$ be the set of tilings in which some ball does not lie
within a supertile of any level $n$. $S$ has zero measure with respect
to any translation invariant measure on $X_\phi$.
\vs.05 \nd
Proof.\ We only give the proof for dimension $m=2$. Note first that
the boundary of a supertile of infinite level must be either a line, or
have a single vertex, since it is tiled by supertiles of all levels and
therefore cannot contain a finite edge. Furthermore, for a given
substitution system there is a constant $K$ such that no tiling in it
contains more than $K$ vertices of supertiles of infinite level;
specifically, one can take $K=2\pi/p$ where $p$ is the smallest angle
of any of the vertices of the tiles.

Next we fix some orthogonal coordinate system in the plane and
decompose $S$ into disjoint subsets as follows. Let $C = [0,1)
\times [0,1)$ be the ``half open'' unit edge square in $\R^2$. Let
$C_{t}$ be the translate of $C$ by the vector $t$.  Let $S'$ be the
subset of $S$ consisting of tilings containing vertices of supertiles of
infinite level. For $x\in S'$ we choose a vertex $V(x)$ by
lexicographic order: we choose that vertex which in the given
coordinate system has the largest first coordinate; if there is more
than one with that coordinate we choose the one with the largest
second coordinate. Then we decompose $S'=\cup_{t\in
\Z^2}S'(t)$, where $x\in S'\cap S'(t)$ if $V(x) \in C_t$. It is easy to
see that each $S'(t)$ is measurable, and that they are translates
of one another so they must have zero measure with respect to any
translation invariant measure.

The tilings $x\in S/S'$ contain two supertiles of infinite level, each
occupying a half plane. Next we decompose $S/S'=\cup_{j,k\in
\Z}\;\sigma_j\cup \sigma'_k$ where $x\in S/S'\cap\sigma_j$ if the
boundary between the supertiles of infinite level crosses the first axis 
in $[j,j+1)$, and $x\in S/S'\cap\sigma'_k$ if the boundary
between the supertiles of infinite level is parallel to the first axis and
crosses the second axis in $[k,k+1)$. Note that all
sets $\sigma_j$ are translates of one another, and all sets
$\sigma'_k$ are translates of one another, so $S/S'$ has zero measure
with respect to any translation invariant measure.\qed
\vs.1 \nd
{\bf Theorem 2.}\ For any substitution tiling system $X_\phi$
satisfying iv), there exists $\epsilon_0>0$ such that
for all $\epsilon \in(0, \epsilon_0)$, and for almost
every tiling $x\in X_{\phi}$, $\O^\epsilon(x)$ is c-equivalent to
$\O_j(x)$ for some (and therefore any) $j$. Up to conjugacy, $\O_j(x)$ is
independent of $x$. Furthermore, if $\phi$ satisfies v) then
$\O^\epsilon(x) = \O_j(x)$ for some (and therefore any) $j$.
\vs.1 \nd
Proof.\ \ From Theorem 1b it follows that, for small $\epsilon$,
$\O^\epsilon(x)$ is contained in $\O_j(x)$, where $j$ is the type of any
of the tiles of $x$ which meet the origin. On the other hand, let $N$
correspond to $\epsilon$ in Theorem 1a. By Lemma 6, for almost
every $x$ there is some $n$ such that the tiles which intersect
$B_{N}$ are contained in some supertile $b$ of level $n$ in $x$. Let
$k$ be the type of $b$. It follows from Theorem 1a that $\O_k(x')$
is a subgroup of $\O^\epsilon(x)$, where $x'=\phi^{-n}x$.  But
$\O_k(x')$ and $\O_j(x)$ are c-equivalent, so $\O_j(x)$ is conjugate to a
subgroup of $\O_k(x')$, and therefore is conjugate to a subgroup of
$\O^\epsilon(x)$.  So $\O_j(x)$ and $\O^\epsilon(x)$ are c-equivalent.
By Lemma 1, $\O_j(x)$ is, up to conjugacy, independent of $x$.
If $\phi$ satisfies v), then $\O_k(x')=\O_k(x)=\O_j(x)$.  Since
$\O_j(x) = \O_k(x') \subset \O^\epsilon(x) \subset \O_j(x)$,
$\O^\epsilon(x)=\O_j(x)$. \qed
\vs.1 \nd
{\bf Corollary 1.} For each substitution tiling system satisfying iv)
the group $\O^\epsilon(x)$ is uniquely defined up to c-equivalence,
for almost all tilings $x$, and all small enough $\epsilon$, thus the
c-equivalence class of the group is an invariant for
equivalence. Furthermore, among substitution tiling
systems also satisfying v), the conjugacy class of this
subgroup of $SO(m)$ is an invariant for equivalence.
\vs.1 \nd
{\bf 3.\ Abstract Substitution Systems}
\vs.1
In going from Lemma 2 to Theorem 2 we see that we can associate with
each substitution tiling system a c-equivalence class of subgroups of
$SO(m)$ in a reasonably fundamental way. We are now ready to relax
the hypotheses.
\vs.1 \nd
{\bf Definition.} A ``substitution (dynamical) system'' is a quadruple
$(X, T, \phi, |\phi|)$ consisting of a compact metric space
$X$ on which there is a continuous action $T\,:\,(g,x)\in
G_E\times X\longrightarrow T^gx\in X$ of $G_E$ and a homeomorphism 
$\phi\,:\,X\to X$ such that
$\phi[T^gx]=T^{\phi(g)}x$ for all $x$, where $\phi(g)$ is the
conjugate of $g$ by the similarity of Euclidean space consisting of
stretching about the origin by $|\phi|>1$.
\vs.1
Substitution tiling systems are special cases of substitution systems.
The map $\phi$ is not intrinsic to the substitution tiling system $(X_\phi,
T)$ since, for tiling systems, $\phi$ and $\phi^k$
lead to the same set of tilings; so equivalence of such systems
should not be required to intertwine the actions of the maps $\phi$. The objects
$S^\epsilon(x)$, $R^\epsilon(x)$ and $\O^\epsilon(x)$ are well defined
in our abstract setting. Motivated by the last section, we use the
following notion of equivalence.
\vs.1 \nd
{\bf Definition.} The substitution systems $(X^1, T^1, \phi^1, |\phi^1|)$ 
and $(X^2, T^2, \phi^2, |\phi^2|)$ are ``equivalent'' if there are 
subsets $Y_j\subset X^j$, invariant under $T^j$ and of
measure zero with respect to all translation invariant Borel
probability measures on $X_{\phi^j}$, 
and a one-to-one, onto, Borel bimeasureable ``isomorphism'' 
$\tau\,:\,X^1\!-\!Y_1\to X^2\!-Y_2\!$,
such that $\tau\circ T^1=T^2\circ \tau$.
Furthermore, $\tau$
must respect the ``local contracting directions''
$S^\epsilon(x)$. Respecting the local contracting directions means
that, for each $x\in X^1\!-\!Y_1$, $\epsilon
>0$ and $\epsilon' >0$, there exist $\tilde \epsilon >0$ and $\tilde
\epsilon' >0$ such that $\tau[S^{\tilde \epsilon}(x)]
\subset S^{\epsilon}(\tau x)$, and $\tau^{-1}[S^{\tilde
\epsilon'}(\tau x)] \subset S^{\epsilon'}(x)$. 
\vs.1
It is easy to see that for the special case of substitution tiling
systems this notion of equivalence reduces to that
previously defined. We will now introduce an invariant for
equivalence which reduces to the class of subgroups of $SO(m)$ we
found for substitution tiling systems.
We note that this allows us to
generalize our discussion of substitution tiling systems to include
tiling systems which do not quite fit the conditions of
\S 2. In particular, our analysis apply to
the various versions of Penrose tilings of the plane, such as the kite
\& dart tilings, both the substitution version and the version with
edge markings [Gar, Ra3], and to the various tilings discussed in
[G-S, Moz].

We will need to introduce a few more definitions. Given two subgroups
$G_1$ and $G_2$ of $SO(m)$ we write $G_1\prec G_2$ if $G_1$ is
conjugate (by an element of $SO(m)$) to a subgroup of $G_2$. The
binary relation $\prec$ lifts in an obvious way to a partial ordering
on the set of c-equivalence classes. We denote by ``lower bound'' to a
set $\tilde S$ of subgroups of $SO(m)$ any c-equivalence class of
groups $G$ each of which satisfies $G\prec S$ for all $S\in \tilde
S$. It is almost immediate that $\O^\epsilon(x)\prec
\O^{\epsilon'}(x)$ if $\epsilon < \epsilon'$.  For each $x\in X$
we define $\G(x)$ as the set of all lower bounds of the family
$\{\O^{\epsilon}(x)\,:\,\epsilon>0\}$; it is nonempty since it
contains $\{e\}$. Note that the set $\G(x)$ is an invariant for
substitution systems -- if $\tau$ is an isomorphism then $\G(\tau
x)=\G(x)$ for almost every $x$. For substitution tiling systems, the
sets $\G(x)$ have unique greatest elements which are constant for
almost every $x$ with respect to every translation invariant 
measure. In the latter case, where $\G(x)$ has an almost everywhere
constant greatest element, we denote this greatest element by
$[\O](X)$. Note that $[\O](X)$ is a c-equivalence class,
unlike $\O^\epsilon(x)$, which is a specific group.  We have thus
generalized the analysis of substitution tiling systems to the more
general setting.

As with substitution tiling systems, we can avoid the use of
c-equivalence classes for systems with a special property.

\nd {\bf Property P:} For almost every $x$ there exists an 
$\epsilon>0$ such that, if $0<\epsilon'<\epsilon$, then
$\O^{\epsilon'}(x)=\O^\epsilon(x)$.

Note that, by Theorem 2, any substitution tiling system that satisfies
v) also satisfies Property P.  If a substitution system
satisfies Property P, we can define $\O(x)$ to be
$\O^\epsilon(x)$ for $\epsilon$ sufficiently small.  If the conjugacy
class of $\O(x)$ is almost everywhere constant, we define $[\O]_0(X)$
to be that conjugacy class.  The previously defined $[\O](X)$ is,
of course, the c-equivalence class of $[\O]_0(X)$.

\nd {\bf Theorem 3.}\ \ Suppose $(X^1,T^1,\phi^1,|\phi^1|)$ and 
$(X^2,T^2,\phi^2,|\phi^2|)$ are
equivalent substitution systems, with the notation of the definition.
Then if $(X^1,T^1,\phi^1,|\phi^1|)$ satisfies Property P so does 
$(X^2,T^2,\phi^2,|\phi^2|)$.
Furthermore, for almost every $x \in X^1$,
$\O(\tau x)=\O(x)$.
In particular, if $\O(x)$ is almost everywhere
constant up to conjugacy then $\O(\tau x)$ is almost everywhere
constant up to conjugacy and $[\O]_0(X^2)=[\O]_0(X^1)$.
\vs.1 \nd
Proof.\ Let $x$ be a generic point of $X^1$. Since $(X^1,T^1,\phi^1,|\phi^1|)$ has
Property P we can find $\epsilon_0>0$ such that, for $0 <
\epsilon < \epsilon_0$, $\O^\epsilon(x) = \O^{\epsilon_0}(x) =
\O(x)$. From the equivalence we can find $\tilde
\epsilon$ such that $\O^{\tilde \epsilon}(\tau x) \subset
\O^{\epsilon_0}(x)$.  Now let $\epsilon_0'=\tilde \epsilon$.  We will
show that, for any $0<\epsilon' \le \epsilon_0'$, $\O^{\epsilon'}(\tau
x) = \O^{\epsilon'_0}(\tau x) = \O(x)$.  From this it will follow that
$(X',\phi')$ has Property P and that $\O(\tau x)=\O(x)$.

Fix any $0< \epsilon' \le \epsilon_0'$.  Since $\tau$ is an isomorphism
there exists $\tilde \epsilon'>0$ such that $\O^{\tilde \epsilon'}(x)
\subset \O^{\epsilon'}(\tau x)$.  But $\O^{\epsilon'}(\tau x) \subset
\O^{\epsilon_0'}(\tau x) \subset \O^{\epsilon_0}(x)$.  If $\O^{\tilde
\epsilon'}(x) = \O^{\epsilon_0}(x)$ then all the inclusions must be
equalities, and we are done.  So it suffices to show  $\O^{\tilde
\epsilon'}(x) = \O^{\epsilon_0}(x)$.  If $\tilde \epsilon' \le
\epsilon_0$ this follows from the definition of $\epsilon_0$.  But if
$\tilde \epsilon' > \epsilon_0$ then $O^{\epsilon_0}(x) \subset
O^{\tilde \epsilon'}(x)$, so  $\O^{\tilde
\epsilon'}(x) = \O^{\epsilon_0}(x)$. \qed
\vs1 \nd
{\bf 4.\ Examples and Analysis of the Invariant} 
\vs.1
For the pinwheel $[\O]_0(X_{\phi_p})$ is the group generated by
rotations by $\pi/2$ and $2\arctan(1/2)$; for the variant of the
pinwheel $[\O]_0(X_{\phi_v})$ is the group generated by rotations by
$\pi/2$ and $2\arctan(1/3)$ [RaS]. It is clear that these are
distinct, so the substitution tiling systems are not equivalent.

Rotations appear in more interesting ways in 3 dimensional tilings,
for example the quaquaversal and dite \& kart substitution tiling
systems, defined in [CoR] and [RaS] respectively. These systems both
satisfy v) and therefore Property P. Let $R_x^\theta$ be
a rotation about the $x$ axis by an angle $\theta$, with similar
notation for other axes.  If we denote by $G(p,q)$ the subgroup of
$SO(3)$ generated by $R_x^{2\pi/p}$ and $R_y^{2\pi/q}$, it can be
shown [CoR, RaS] that $[\O]_0(X_{\phi_q})$ is the conjugacy class of
$G(6,4)$ for the quaquaversal tilings and $[\O]_0(X_{\phi_{d\&k}})$ is
the conjugacy class of $G(10,4)$ for the dite \& kart tilings.  We
shall see that $G(6,4)$ and $G(10,4)$ are not conjugate (indeed not
even c-equivalent) by using the following obvious fact: if the groups
$G$ and $G'$ are conjugate (or c-equivalent) and one of them has an
element of order $m$ (finite or infinite) then the other must have an
element of order $m$.

\vs.1 
\nd {\bf Structure Theorem for G(p,q)}\ [RaS]

\nd $a)$ If $p,q\ge 3$ are odd, then $G(p,q)$ is isomorphic to
the free product
$$ \Z_p * \Z_q =\ <\a,\b\,:\,\a^p,\,\b^q>. \eqno 5)$$

\nd $b)$ If $p\ge 4$ is even and $q\ge 3$ is odd, then $G(p,q)$ has the
presentation 
$$<\a,\b\,:\,\a^p,\,\b^q,\,(\a^{p/2}\b)^2>. \eqno 6)$$

\nd $b)$ If $p\ge 4$ is even and $q=2s$, $s\ge 3$ odd, then $G(p,q)$ 
has the
presentation 
$$<\a,\b\,:\,\a^p,\,\b^q,\,(\a^{p/2}\b)^2,\,(\a\b^s)^2>. 
\eqno 7)$$

\nd $d)$ If $4$ divides both $p$ and $q$, then $G(p,q)=G([p,q],4)$,
where $[p,q]$ denotes the least common multiple of $p$ and $q$.

\nd $e)$ If $4$ divides $m$, then $G(m,4)$ has the presentation
$$<\a,\b,:\,\a^m,\,\b^4,\,(\a^{m/2}\b)^2,\,(\a\b^2)^2,
(\a^{m/4}b)^3>. \eqno 8) 
$$
In cases a), b) and c), the isomorphism between the abstract
presentation and $G(p,q)$ is given by $\a \mapsto R_x^{2\pi/p}$, $\b
\mapsto R_y^{2\pi/q}$.  In case e) the isomorphism is similar.
\vs.1 \nd
{\bf Theorem 4.}\ \ a) If 4 does not divide both $p$ and $q$ then the
orders of elements of finite order in $G(p,q)$ are $\{\hbox{factors of
}p\}\cup \{\hbox{factors of }q\}$; b) If 4 divides both $p$ and $q$
then the orders of elements of finite order in $G(p,q)$ are
$\{\hbox{factors of }[p,q]\}\cup \{3\}$.
\vs.1 \nd
{\bf Corollary 2.}\ \ If $p$ and $q$ are not both divisible by $4$,
and $p'$ is not a factor of $p$ or $q$, then $G(p,q)$ and $G(p',q')$
are not c-equivalent.
\vs.1 \nd
{\bf Corollary 3.}\ \ The quaquaversal and dite \& kart systems are
not equivalent.

\nd
Proof of the theorem.\vs0 \nd {\bf a)} Assume $g\in G(p,q)$ has finite
order $\ne 1$. We know $g$ can be expressed in one of the forms
$g=A^{a_1}B^{b_1}A^{a_2}\cdots A^{a_{n+1}}$,
$g=B^{b_1}A^{a_1}B^{b_2}\cdots B^{b_{n+1}}$,
$g=B^{b_1}A^{a_1}B^{b_2}\cdots A^{a_n}$ or
$g=A^{a_1}B^{b_1}A^{a_2}\cdots B^{b_n}$, with all $0< a_j < p$,\ $0 <
b_j < q$ and $n\ge 1$.  Within the class of $g'\in G(p,q)$ which are
conjugate to $g$ with respect to $G(p,q)$, we assume that $n$ is
minimal. Assume $n\ge 2$; we will obtain a contradiction.  Since one
could conjugate with $A^{a_1}$, the form $A^{a_1}B^{b_1}A^{a_2}\cdots
A^{a_n}$ can be exchanged for $B^{b_1}A^{a_2}\cdots A^{a_n}$ and since
one could conjugate with $B^{b_1}$, the form
$B^{b_1}A^{a_1}B^{b_2}\cdots B^{b_n}$ can be exchanged for
$A^{a_1}B^{b_2}\cdots B^{b_n}$, and the form
$B^{b_1}A^{a_1}B^{b_2}\cdots A^{a_n}$ for $A^{a_1}B^{b_1}A^{a_2}\cdots
B^{b_n}$.  So we assume that the form is $A^{a_1}B^{b_1}A^{a_2}\cdots
B^{b_n}$.  If $a_j=p/2$ (resp.\ $b_j=q/2$) for any $j$ then we could
use the relation $A^{p/2}B^{b}=B^{-b}A^{p/2}$ (resp.\
$B^{q/2}A^a=A^{-a}B^{q/2}$) to reduce the value of $n$; thus these
values of $a_j$ (or $b_j$) cannot occur. But then, by the structure
theorem, $g$ has infinite order, which is a contradiction. Thus $n$
must equal 1, and $g$ can be assumed to be of the form $A^{a_1}$,
$B^{b_1}$ or $A^{a_1}B^{b_1}$. Considering
$A^{a_1}B^{b_1}A^{a_1}B^{b_1}\cdots A^{a_1}B^{b_1}$, the only way
$A^{a_1}B^{b_1}$ could have finite order is if $a_1=p/2$ or $b_1=q/2$,
in which case $g$ has order 2, and 2 is a factor of $p$ or
$q$. Finally, the elements $A^{a_1}$ can have as orders any factor of
$p$ and the elements $B^{b_1}$ can have as orders any factor of $q$.
\vs.1
\nd {\bf b)} If $p$ and $q$ are divisible by 4 then $G(p,q)=G([p,q],4)$, so
we consider $G(m,4)$ with $m$ divisible by 4. Using the presentation
8), we can put any 
$g\in G(m,4)$ in the form $WST^{a_1}ST^{a_2}\cdots
ST^{a_n}E$ with $S=R_y^{2\pi/4}$, $T=R_x^{2\pi/m}$, $n\ge 0$, $a_j\ne
km/4$ and with both $W$ and $E$ in the cube group $G(4,4)$. Assume $g$ has
finite order $\ne 1$ and that in its conjugacy class (which of course
all have the same order), the smallest value of $n$ in the above
representation is $\ge 2$. (We will obtain a contradiction to this.)
By conjugation we eliminate $W$ from $g$.

Now $G(4,4)$ can be partitioned: $G(4,4)=H_1\cup H_1S\cup H_1SU$, where
$U=R_x^{2\pi/4}$ and 
$H_1$ is the 8 element subgroup generated by $S^2$ and $U$. In detail,
$$H_1=\{1,U,U^2,U^3,S^2,S^2U,S^2U^2,S^2U^3\}. \eqno 9)$$
\nd Some power of $g$ equals the identity element:
$$(ST^{a_1}ST^{a_2}\cdots ST^{a_n}E)(ST^{a_1}ST^{a_2}\cdots
ST^{a_n}E)(\cdots) =e. 
\eqno 10)$$ 
\nd We consider the three cases: i) $E\in H_1$; ii)
$E\in H_1SU$; iii) $E\in H_1S$. 

\nd {\bf i)} The factor $E$ in 10) 
is of the form $(S^2)^aU^b$ with $a=0,1$ and $b=0,1,2,3.$ We alter
10) to
$$[ST^{a_1}ST^{a_2}\cdots ST^{a_n}U^{(-1)^ab}(S^2)^a][ST^{a_1}ST^{a_2}\cdots
ST^{a_n}U^{(-1)^ab}(S^2)^a][\cdots ] =e, \eqno 11)$$ 
\nd or
$$\eqalign{[ST^{a_1}ST^{a_2}\cdots
ST^{[a_n+(-1)^abm/4]}&(S^{2a})][ST^{a_1}ST^{a_2}\cdots \cr
&\cdots ST^{[a_n+(-1)^abm/4]}(S^{2a})]\cdots
=e,\cr}\eqno 12)$$
\nd and we know [RaS] this cannot be the case. So we cannot have $E\in H_1$.

\nd {\bf ii)} $E$ is now of the form $(S^2)^aU^bSU$ with 
$a=0,1$ and $b=0,1,2,3.$ We now alter 10) to
$$ST^{a_1}ST^{a_2}\cdots 
ST^{[a_n+(-1)^abm/4]}(S^2)^aSUST^{a_1}ST^{a_2}\cdots
ST^{a_n}E\cdots =e.\eqno 13)$$
\nd Using $SUS=USU$, 13) becomes
$$\eqalign{ST^{a_1}S&T^{a_2}\cdots
ST^{[a_n+(-1)^a(b+1)m/4]}S^{2a}ST^{a_1+m/4}ST^{a_2}\cdots \cr
&\cdots ST^{[a_n+(-1)^a(b+1)m/4]}S^{2a}SU\cdots =e.}\eqno 14)$$

\nd Again, we know [RaS] 
this cannot be the case. So we cannot have $E\in H_1SU$.

\nd {\bf iii)} We cannot have $E\in H_1S$ and $n \ge 2$. For if we
represent conjugacy by $\cong$,
$$\eqalign{g\cong ST^{a_1}ST^{a_2}\cdots ST^{a_n}E &
=ST^{a_1}ST^{a_2}\cdots ST^{a_n}(S^2)^aU^bS \cr & \cong
T^{a_1}ST^{a_2}\cdots ST^{[a_n+(-1)^abm/4]}(S^2)^{a+1} \cr & \cong
ST^{a_2}\cdots ST^{[a_n+(-1)^abm/4+(-1)^{a+1}a_1]}(S^2)^{a+1}, \cr}
\eqno 15)$$
\nd and $g$ is conjugate to a word with smaller $n$.

Thus $n=0$ or $n=1$.  $n=0$ means $g\in G(4,4)$, and these have orders
$1,2,3,4$. $n=1$ means $g$ is of the form $ST^{a_1}E$ where $a_1\ne
km/4$ and $E\in G(4,4)$. We again consider the three cosets to which
$E$ may belong. As before we see that cases i) and ii) lead to
infinite order for $g$. But in case iii) $g$ is conjugate to
$T^{a_1}(S^2)^aU^bS^2\cong T^c(S^2)^d$, which can have as orders the
factors of $m$. \qed


\nd {\bf 5.\ Conclusion}

We have been concerned with substitution tilings of Euclidean spaces,
and have defined an invariant for them related to the group
generated by the relative orientations of the tiles in a tiling. This
feature is captured in an intrinsic way by means of a contractive
behavior of the substitution. It is unrelated to other features of
tiling systems, such as their topology, and we introduce the notion of
substitution dynamical system to emphasize the features associated
with the invariant.

To distinguish examples, for instance to distinguish the quaquaversal
tilings from the dite \& kart tilings, requires consideration of
2-generator subgroups of $SO(3)$, in particular the orders of elements
of such subgroups, which we analyze.

\nd {\bf Acknowledgements.} We are pleased to thank Ian Putnam for
useful discussions.


\vfill\eject
\centerline{{\bf References}}
\vs.2 \nd
[AnP]\ J.E.\  Anderson and I.F.\ Putnam, Topological invariants for
substitution tilings and their associated $C^\ast$-algebras,
{\it Erg.\ Th.\ Dyn.\ Sys.}, to appear.
\vs.1 \nd
[CoR]\ J.H.\ Conway and C.\ Radin, Quaquaversal tilings and rotations,
Inventiones math, to appear. Available
from the electronic preprint archive  
\nd \vskip -6 pt \hs -.33 
mp\_arc@math.utexas.edu as 95-425.
\vs.1 \nd
[CoK]\ E.\ Coven and M.\ Keane, The structure of substitution minimal sets,
{\it Trans.\ Amer.\ Math.\ Soc.} 162 (1971), 89-102.
\vs.1 \nd
[DwS]\ S.\ Dworkin and J.-I Shieh, Deceptions in quasicrystal growth,
{\it Commun.\ Math.\ Phys.} 168 (1995), 337-352.
\vs.1 \nd
[Gar]\ M.\ Gardner, Extraordinary nonperiodic tiling that enriches the theory
of tiles, {\it Sci.\  Amer.\ } January 1977, 116-119.
\vs.1 \nd
[G-S]\ C.\  Goodman-Strauss, Matching rules and substitution tilings,
preprint, University of Arkansas, 1996.
\vs.1 \nd
[GrS] B.\ Gr\"unbaum and G.C.\ Shephard, {\it Tilings and Patterns},
Freeman, New York, 1986.
\vs.1 \nd
[Kel]\ J.\ Kellendonk, Non-commutative geometry of tiligs and gap 
labelling, {\it Rev.\  Mod.\ Phys.}, to appear.
\vs.1 \nd
[Ken]\ R.\ Kenyon, Inflationary tilings with a similarity structure, 
{\it Comment.\ Math.\ Helv.} 69 (1994), 169-198. 
\vs.1 \nd
[LaW]\ J.C.\ Lagarias and Y.\ Wang, Self-affine tiles in $\R^n$,
{\it Adv. Math.}, to appear.
\vs.1 \nd
[Min]\ J.\ Mingo, $C^\ast$-algebras associated with one dimensional almost 
periodic tilings, {\it Commun.\ Math.\ Phys}, to appear.
\vs.1 \nd
[Moz]\ S.\ Mozes, Tilings, substitution systems and dynamical systems
generated by them, {\it J. d'Analyse Math.}\ 53 (1989), 139-186.
\vs.1 \nd
[CEP]\ H.\ Cohn, N.\ Elkies and J.\ Propp,
Local statistics for random domino tilings of the Aztec diamond,
{\it Duke Math.\  J.}, to appear.
\vs.1 \nd
[Ra1] C.\ Radin, The pinwheel tilings of the plane, {\it Annals of Math.}
139 (1994), 661-702.
\vs1 \nd
[Ra2]\ C.\ Radin, Space tilings and substitutions, {\it Geometriae Dedicata}
55 (1995), 257-264.
\vs.1 \nd
[Ra3]\ C.\ Radin, Miles of Tiles, Ergodic theory of ${\Z}^d$-actions,
{\it London Math. Soc. Lecture Notes Ser.\ } 228 (1996), Cambridge Univ.\ 
Press, 237-258.
\vs.1 \nd
[RaS]\ C.\ Radin and L.\ Sadun, Subgroups of SO(3) associated with tilings,
{\it J. Algebra}, to appear. Available
from the electronic preprint archive mp\_arc@math.utexas.edu as 96-4.
\vs.1 \nd
[RaW]\ C.\ Radin and M.\ Wolff, Space tilings and local isomorphism,
{\it Geometriae Dedicata} {42} (1992), 355-360.
\vs.1 \nd
[Rob]\ E.A.\ Robinson, The dynamical properties of Penrose tilings,
{\it Trans.\ Amer.\ Math.\ Soc.}, to appear.
\vs.1 \nd
[Sad]\ L.\ Sadun, Some generalizations of the pinwheel tiling,
{\it Disc.\ \& Comp.\ Geo.}, to appear. Available
from the electronic preprint archive mp\_arc@math.utexas.edu as 96-97.
\vs.1 \nd
[Sch]\ K.\ Schmidt, Tilings, fundamental cocyles and fundamental groups
of symbolic $\Z^d$-actions, preprint, Erwin Schr\" odinger Institute, 1996
\vs.1 \nd
[Sen]\ M.\ Senechal, {\it Quasicrystals and geometry}, Cambridge
University Press, Cambridge, 1995.
\vs.1 \nd
[Sol]\ B.\ Solomyak, Non-periodic self-similar tilings have the unique 
composition property, preprint, University of Washington, 1996.
\vs.1 \nd
[Tha]\ L.\ Thang, Local rules for pentagonal quasi-crystals,
{\it Disc.\ \& Comp.\  Geo.}, 14 (1995), 31-70.
\vs.1 \nd
[Wa1]\ H.\ Wang, Proving theorems by pattern recognition II, {\it Bell 
Systs.\ Tech.\ J.}\ {40} (1961), 1-41.
\vs.1 \nd
[Wa2]\ H.\ Wang, Games, logic and computers, {\it Sci.\ Am.\ (USA)}
(November 1965), 98-106.

\end